\documentclass[12pt,leqno]{article}
\usepackage{graphicx}
\begin{document}

\begin{center}
{\bf
Srishti Dhar Chatterji (1935-2017)
\vskip 12pt
In Memoriam
\vskip 12pt
V.S.~Varadarajan and Robert C.~Dalang
}
\vskip 12pt
University of California Los Angeles\\ and\\ \'Ecole Polytechnique F\'ed\'erale de Lausanne

\end{center}
\vskip 16pt

\hfill {\small J\"ons: {\em  Feel, to the very end, the triumph of being alive}!

\hfill Karin:  {\em Quiet, quiet}!

\hfill J\"ons:  {\em Yes, but under protest.}
\vskip 12pt

\hfill Ingmar Bergman, The Seventh Seal.}
\vskip 16pt

Professor Srishti Dhar Chatterji passed away on September 28, 2017, in Lausanne, Switzerland, most suddenly and unexpectedly, after a very brief illness.

\section{Youth and early career}

S.D. Chatterji, known to most as ``Chatterji," or more simply, ``Chat", was born on June 29, 1935, in Parbatipur near the township of Sheakhala, approximately 40 km to the north-west of Calcutta, India. This was the hometown of his mother Shivani Chatterjee. During Chat's first few years, the family lived in Sheakhala, the hometown of his father Kali Pada Chatterjee (who lived up to age 98), where his ancestors had resided for about 300 years, and Chatterji's grandfather was a notable person there, of the Brahmin caste. The name Chatterji (as it was spelled on his first passport), or Chatterjee, is a clan name, associated to a large group of people from Bengal. It is an anglicization of the Indian name Chattopadhyaya (``upadhyaya" means teacher, and is an honorific).

In 1939, for both economic and personal reasons, the family moved to Lucknow, the main city in the state of Uttar Pradesh, about 1000 km away. Chat had five younger brothers and one younger sister, several of which had successful careers of their own, but of whom only the brothers Shyam Sunder Chatterjee and Dharani Dhar Chatterjee outlive him. In Lucknow, Chat received tough home-schooling from his father, and listened often to the BBC radio (which may have given him his excellent English accent). The testimonial article of F.~Labourie in this volume describes Chatterji's further schooling.

Chatterji's mother tongue was Bengali, but he later learned, and was fluent in, Hindi, Urdu, English, Italian, German, French and Danish. He also wrote fluently in these languages, as can be seen from his list of publications.

  He entered the University of Lucknow in 1950 at the young age of 15, where he obtained a Bachelor of Science in Physics, Mathematics and Statistics in 1952, and a Master degree in Mathematical Statistics two years later. As a youth, Chatterji was also accomplished in sports, particularly cricket and squash. He even considered becoming a professional cricket player.
	
	In 1954, Chatterji entered the Indian Statistical Institute (ISI) in Calcutta, India, as a graduate student, where he met the famous U.S. mathematician Norbert Wiener, who was a visiting professor there during 1955-56. Wiener had the young Chat read papers of the famous probabilist Kai-Lai Chung, gave him some research projects and encouraged him to go and study in the United States. In his second year at the ISI, Chat met the entering graduate student Veeravalli S.~Varadarajan, and the two entered into a strong friendship that continued to the end. As Varadarajan recounts, ``together with R. Pakshirajan, we formed a trio, trying to understand Doob's great book on Stochastic Processes \cite{doob}\footnote{Papers in Chatterji's list of publications are labelled with numbers; other references are labelled with letters.}, without any supervisory guidance. To me, then a green lad of 19, Chat was the epitome of knowledge and sophistication, and that attitude to him has always been with me, even after so many years. Chat was impatient at our modest progress." In summer of 1956, Chat left for the U.S.A., initially to work with Kai Lai Chung at Syracuse University in upstate New York: see the testimonial article of F.~Labourie for details of this voyage.
	
	After two years at Syracuse, Chatterji moved to Michigan State University, where he completed his thesis in 1960 under Professor Charles Kraft. He mentioned much later to Varadarajan that the late Professor Gopinath Kallianpur, a famous Indian probabilist, was of much help to him in his thesis on Banach space-valued martingales \cite{rd19a}, which was a major work (more about this later).
	
	From 1960 to 1970, Chatterji worked at several universities in various parts of the world. His first position (1960-62) was as Lecturer in the Department of Mathematics at the University of New South Wales, Australia. He then returned as an Assistant Professor to the Department of Statistics in Michigan State University (1962-63), before spending a year in the Mathematics Research Center of the University of Wisconsin, Madison. Chatterji then moved to Europe, first as a Gastprofessor in the Mathematisches Institut der Universit\"at Heidelberg, Germany (1964-66).
	
	He was then hired as Lektor in the Statistisk Institut at the University of Copenhagen, Denmark (1966-67). Chat's next position was in the Advanced Studies Center, Institut Battelle, a private research center in Geneva, Switzerland (1967-68), where his future colleague Renzo Cairoli (1931-1994) also worked. He then returned for a year to Canada, as a Professeur invit\'e in the D\'epartement de math\'ematiques, Universit\'e de Montr\'eal (1968-69), before going back as a Lektor to the Matematisk Institut of the University of Copenhagen (1969-70).
	
	Finally, his last professional move was to Switzerland in 1970, where he took up a position as Professeur ordinaire (formally in statistics) in the D\'epartment de math\'ematiques of the \'Ecole Polytechnique F\'ed\'erale de Lausanne (EPFL), and remained there until his retirement in 2000. This move was not motivated by financial considerations, since Chat's salary in Copenhagen was higher than his initial salary at EPFL, but rather by his common research interests in probability theory with R.~Cairoli, who had also just moved to EPFL, and by personal considerations.
\begin{figure}
\begin{center}
\includegraphics[height=7.5cm]{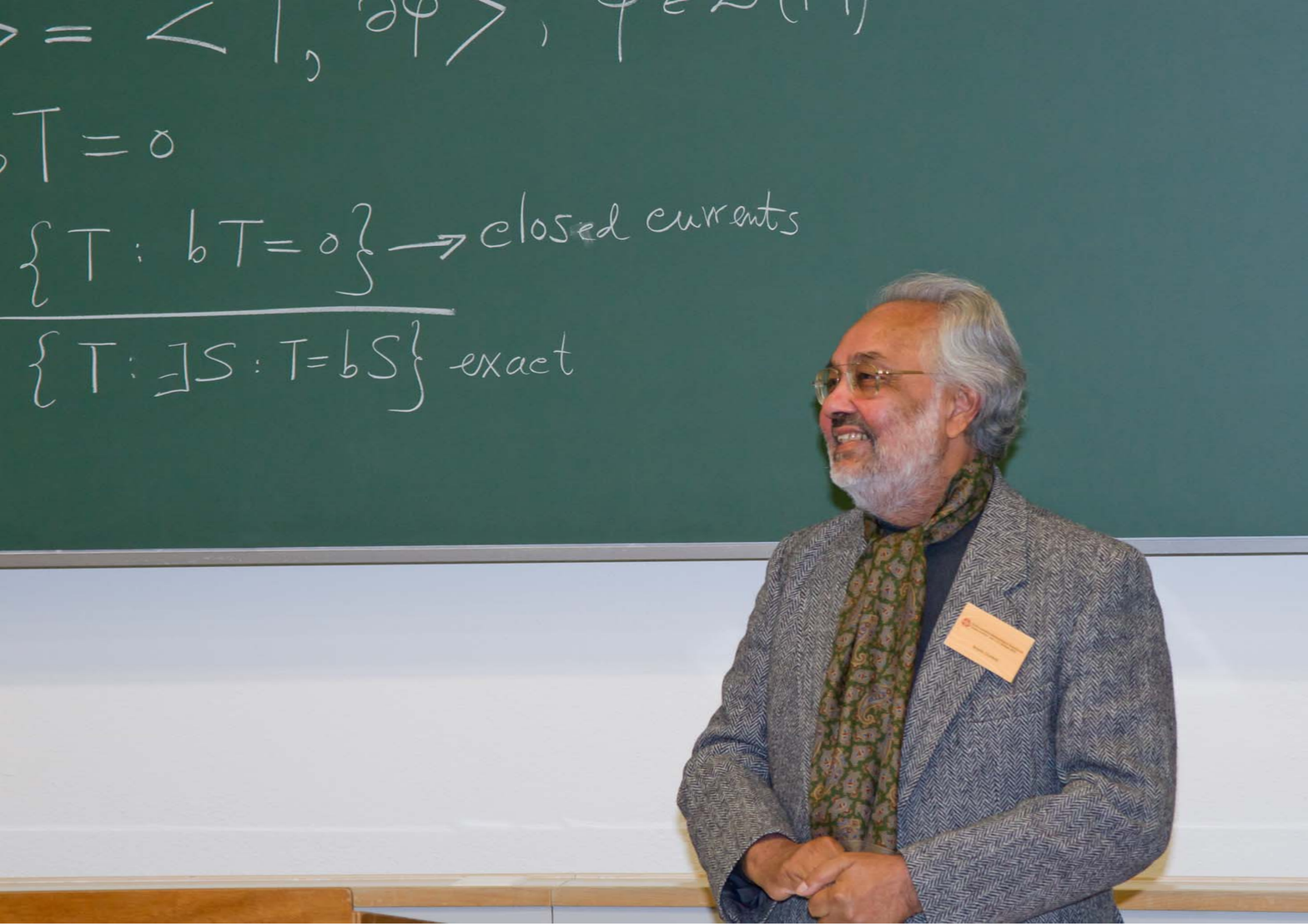}
\end{center}
\caption{Chatterji in 2010, during the centennial celebrations of the Swiss Mathematical Society in Berne, Switzerland. Photo by Prof.~Burchard Kaup (Universit\'e de Fribourg), printed with his permission.}
\end{figure}

\section{Mathematical work}
Chatterji was an expert who worked in measure theory, probability theory, real analysis and functional analysis.
His reputation rests firmly on two foundational theorems that he established. The first concerns convergence theorems for martingales with values in a Banach space $B$, and establishes precisely for which $B$ these convergence theorems are true.
The second is a series of theorems which establish a property called the {\it subsequence principle}: an informal statement of this principle is that given an arbitrary sequence of random variables (satisfying a tightness condition), there exists a subsequence which has all the properties of a sequence of independent and identically distributed (i.i.d.) random variables. It is therefore not possible, by statistical means to distinguish this subsequence from a sequence of i.i.d.~random variables.

	Regarding the first topic, already in his Ph.D.~thesis \cite{rd19a}, Chatterji had started to study the problem of convergence of martingales, for random variables taking values in a Banach space. For real-valued random variables $X_n$ ($n=1,2,\dots$), one says that the sequence $\{X_n\}$ is a martingale if
$$
E|X_n|<\infty
$$
for all $n$, and
$$
E(X_n|X_i,\dots ,X_1) =X_i, \qquad \mbox{for all } i<n.
$$
This concept dates back to gambling casinos in the eighteenth and nineteenth centuries. Heuristically, martingales are the idealization of the notion of \lq\lq fair game.\rq\rq In modern times, it was the U.S. mathematician Joseph L.~Doob who realized their fundamental nature and proved the remarkable theorem that martingales always converge to a limiting random variable almost everywhere or in $L^1$, provided they are uniformly integrable in a natural sense. In more recent times, D.W.~Stroock and S.R.S.~Varadhan \cite{stroockvarad} have shown how the whole theory of diffusion processes can be developed with a remarkable clarity and architectural beauty using martingale theory as a foundation. Doob had asked one of his students, S.C.~Moy, to investigate the case of Banach space-valued random variables. Chatterji took up this area and started the study of the convergence of martingales with values in a Banach space $B$. In his thesis, he proved the convergence theorem when $B$ is reflexive, and constructed counter-examples showing that the result is  not true for arbitrary Banach spaces. He recognized that it should be possible to identify what property of a Banach space $B$ would ensure that $B$-valued martingales converge. Of course, many people were working in this area, but in a definitive paper \cite{c15}, published in 1968, he settled the issue in a decisive manner. Roughly speaking, almost any type of convergence of all $B$-valued martingales is equivalent to the fact that $B$ must have the so-called Radon-Nikodym (RN) property, that is, the Radon-Nikodym theorem is true for $B$-valued countably additive set functions. The RN property is the following: if $\mu$ is a probability measure on the space where the random variables are defined, and $\nu$ is a $B$-valued countably additive set function of finite total variation, then for $\nu$ to be absolutely continuous with respect to $\mu$, namely that $\mu(E)=0\Longrightarrow \nu(E)=0$, it is necessary and sufficient that
$$
\nu (E)=\int _E f\, d\mu
$$
for all measurable sets $E$, where $f:=d\nu/d\mu$ is an integrable $B$-valued function; $f$ is called the {\it Radon-Nikodym derivative} of $\nu$ with respect to $\mu$ and one writes $f=d\nu/d\mu$.  Of course, Chatterji constructed counter examples to illuminate all these equivalences.

Turning to the second topic, inspired (as early as 1970) by results of Steinhaus and Komlos \cite{K} on the strong law of large numbers, Chatterji first formulated the subsequence principle in 1972 \cite{c22}, and carried on the investigation with remarkable power and originality. The subsequence principle states that if a certain quantitative asymptotic property $\Pi$ is valid for any sequence of independent, identically distributed random variables belonging to some integrability class determined by the finiteness of a norm $||\cdot||_L$, an analogous property $\tilde\Pi$ will be valid for a suitable subsequence $\{F_0\}$ of {\it any} sequence $F$ of functions on any probability space such that $\sup \{||f||_L : f\in F\} < \infty$. Moreover, the subsequence can be chosen in such a way that any further subsequence will have the same property $\tilde\Pi$.

	In a series of papers \cite{c26,c27,c28}, all published in 1974, he proved this principle in the case where $\Pi$ is the strong law of large numbers in its many forms, the central limit theorem, and the law of the iterated logarithm.  The point is that the validity of these famous laws is now seen to be true for many sequences which are far from being independent. Chatterji's principle was so general that it was even unclear how to formulate it mathematically in a rigorous manner. This difficult issue was settled in 1977 in the Ph.D.~thesis of David Aldous \cite{aldous}.
	

These two sets of theorems and their consequences, elaborated by him in many papers, display his originality and technical power very clearly. He was very fond of these results, and he enjoyed discussing them in detail. In addition to these two major themes, Chatterji wrote many short papers on delicate aspects of real analysis and functional analysis, often giving very simple proofs of classical theorems whose proofs were very involved.

\section{Mathematical historian, editorial work, and expositor}

As he grew older, Chatterji continued to work in the theory of measure and integration, particularly on differentiation of measures in infinite dimensions \cite{c32,c35,c36} and generalized Riemann integrals for functions with values in a Banach space \cite{c69}. He also found a beautiful and simple proof of the Hausdorff-Young inequality for all locally compact abelian groups \cite{c68}. However, over the years, he slowly gravitated towards the history of mathematics, with special reference to his beloved subjects. He was very interested in and wrote about the work of Albert Einstein \cite{c45}. His three papers on the life and work of Norbert Wiener \cite{c65,c67,c68} are expertly done. He edited, along with H. Wefelscheid, the Selected Papers of G.C.~Young and W.H.~Young \cite{rd15}. He wrote several detailed commentaries for three volumes of the Complete works of Felix Hausdorff \cite{rd16,rd17,rd19}, whose modern version, called {\it The Hausdorff Edition}, was published in 2006. He also edited, and translated into French, a book of pedagogical letters written by Leonhard Euler to a young German princess \cite{rd18}.  He regularly contributed to the Jahrbuch \"Uberblicke Mathematik and was co-editor of 13 volumes of this journal  \cite{rd4}--\cite{rd13}. In 2013, he and Professor Manuel Ojanguren wrote a survey of the work of the famous Swiss mathematician Georges de Rham and its impact on mathematics, Swiss as well as world-wide \cite{c73b} (also reprinted in \cite{c73a}).

   In addition, Chatterji wrote about the work of P.R.~Masani \cite{rd22a,u5}, a mathematician of Indian origin, who had collaborated in a major way with Wiener on multivariate prediction theory. And he has a long discussion on the Borel-Lebesgue exchanges about the origins and priorities of Lebesgue measure theory \cite{u6}).

\section{Chatterji at EPFL}
	
	During his first year in Heidelberg, Chatterji had met his future wife Carla Bolognini.  Carla and Chat married in Copenhagen in March 1967, despite some resistance from both sets of parents. Carla was Swiss, from the Italian-speaking Canton of Ticino, and, naturally, preferred moving back to Switzerland. The couple soon moved into a house in downtown Lausanne, which, at that time, was within walking distance of the EPFL Mathematics Department. Chat became a Swiss citizen in 1982, after going through a strict naturalization process.
	
	At EPFL, Chatterji joined several other mathematicians that he had met at the Institut Batelle, including Michel Andr\'e (1936-2009) and Bruno Zwahlen (1934-2018), as well as his colleague and fellow probabilist Renzo Cairoli, who also became a good friend of his. He developed his research program in probability theory and analysis. Chatterji also participated activily in teaching and in service to the department, as we shall see below. He attracted other high level researchers, including Antonio Gualtierotti (1977 to 1978), who also worked in probability theory, Charles-\'Edouard Pfister (1981 to 1996), a distinguished mathematical physicist, and Jacques Sesiano (1982 to 2000), who was an expert in the history of mathematics.
		
	Chatterji continued to travel widely, and many mathematicians also visited Cairoli and Chatterji at EPFL. John B.~Walsh, who worked with Cairoli, recalls that during social dinners with Chatterji, once all were served with food and wine, they would enjoy the warm conversation of Chatterji and good cheer would flow from Chatterji down the table to the other diners. Varadarajan visited EPFL in 1975 at Chatterji's invitation, to give lectures in the Troisi\`eme Cycle de la Suisse Romande, and came again as a visiting professor in 1994. At one point, Varadarajan had to sign for a railway pass that would enable him to travel freely in all of Switzerland, and there was a question about which class. Varadarajan recounts ``I hesitatingly looked at him, and he replied firmly `professors always travel first class!' He took care of my wife and myself always with great consideration and courtesy. He treated all visitors to EPFL in this way."
	
	Chat and Carla's only child was a daughter, Indira.  Chat was overjoyed when Indira chose to become a professional mathematician, working in Geometric Group Theory and now a professor in the Universit\'e de Nice. He told Varadarajan, perhaps facetiously, that Indira's ideas and professional habits were of the new generation which he did not comprehend, and so he seldom discussed mathematics with her!
\begin{figure}[h]
\begin{center}
\includegraphics[height=7.5cm]{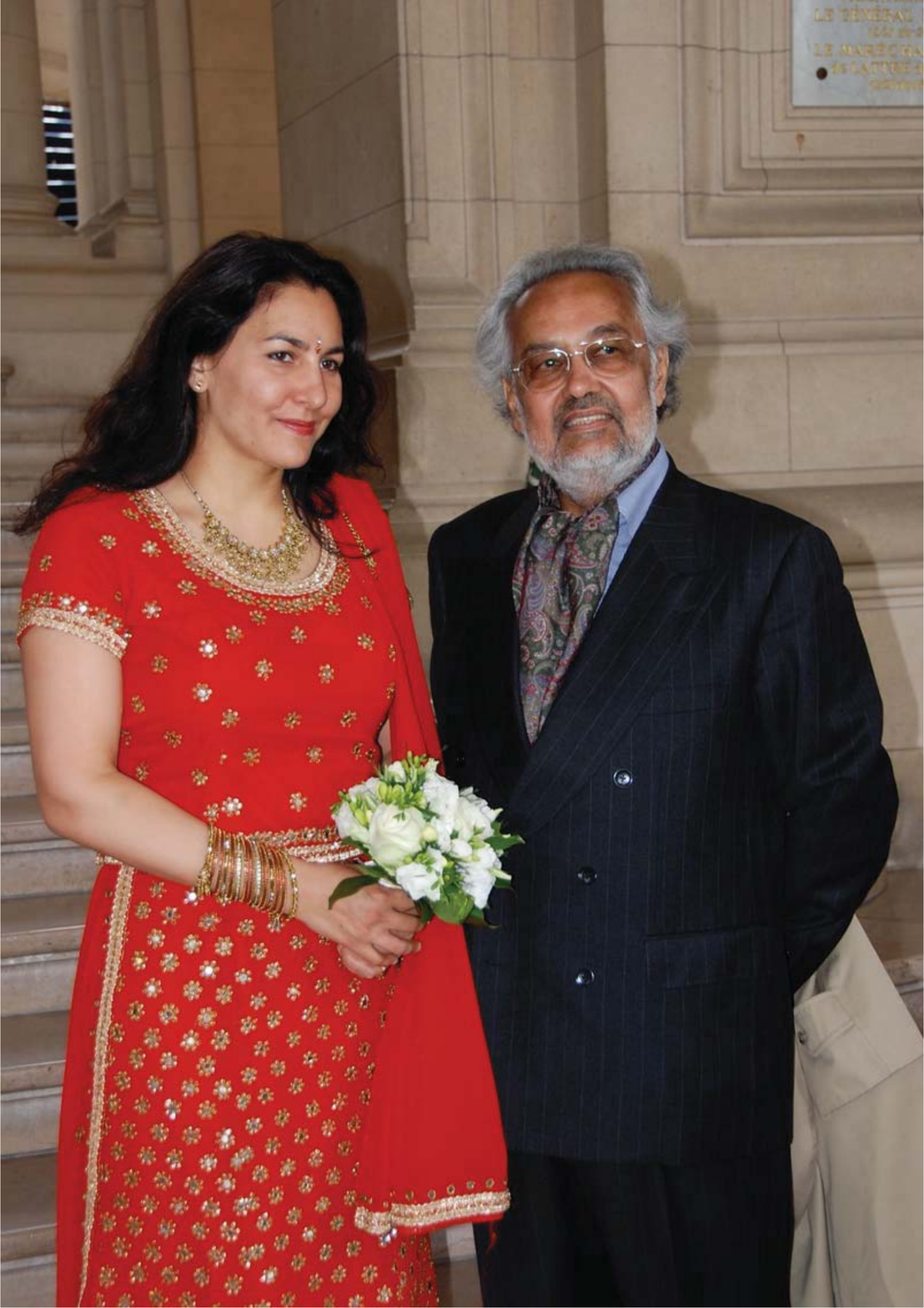}
\end{center}
\caption{Chatterji and his daughter Indira, at her wedding in 2008. Printed with the permission of Indira Chatterji.}
\end{figure}

\section{Teacher}
Chatterji's high level course in the  1971 \'Ecole d'Et\'e de Probabilit\'es de Saint-Flour \cite{c24} was widely recognized for its clarity and elegance, and this contributed greatly to Chatterji's growing reputation. It influenced other leading researchers such as the French mathematician Michel Emery who also appreciated Chatterji's talents as an expositor.

At EPFL, Chatterji taught a ``cours de troisi\`eme cycle" on quantum mechanics in 1973-74, and a course on stochastic processes in the Troisi\`eme Cycle Romand de Math\'ematiques in 1975. These were high level courses for graduate students and other researchers in mathematics.

At the undergraduate level, Chatterji taught many different courses, including ``M\'ethodes math\'ematiques de la Physique (1979-1983)," ``Analyse complexe," ``Analyse harmonique," ``Equations diff\'erentielles," ``M\'ecanique quantique," and, from 1975 to 1983, the first four semesters of the Real Analysis sequence for students in Civil and Rural Engineering. From 1983 to 2000, he taught the second year Analysis course for physicists (also taken by mathematic majors from 1983 to 1986), which was highly appreciated by the students. This led him to write a monumental set of three volumes (in French), totaling over 1'800 pages, titled {\it Cours d'analyse} \cite{b1,b2,b3}, the first of which is dedicated to his father and mother, the second to his wife Carla and daughter Indira, and the third to ``tous mes \'etudiantes et \'etudiants." These books, published by the Presses Polytechniques et Universitaires Romandes between 1996 and 1999, form a valuable high-level reference to subject matters such as multivariate calculus, implicit function theorem, vector calculus, complex analysis, ordinary differential equations, Fourier series and Fourier and Laplace transforms, and partial differential equations.

Chatterji supervised six Ph.D.~students, the list of which is given at the end of this article. His fourth Ph.D.~student, Francesco Russo, is now a professor of mathematics at ENSTA-ParisTech in Paris, France. Like all EPFL mathematics professors, Chatterji supervised numerous Master projects (which, at the time, were called ``projets de dipl\^ome") and semester projects.

As a contribution to continuing education for high school teachers of mathematics, Chatterji taught two courses for this type of audience, one in 1983 for teachers in Suisse Romande, and one in 1986 for teachers in the Canton of Ticino.

\section{Expositiones Mathematicae}
Varadarajan recalls that during the 1970's, Chatterji had several times expressed his feeling that there were not enough journals in the western world that carried expository articles of important areas of mathematics, like the Uspekhi of the (then) Soviet Union. The result of Chatterji's thought processes on this issue was that he decided to rectify this gap himself. He founded {\it Expositiones Mathematicae}, a journal dedicated mainly to expository articles in mathematics, initially owned and published by the Bibliographisches Institut in Mannheim, Germany (and later sold several times to larger and larger companies). The first issue appeared in 1983. He was its founder and Managing Editor from the beginning through 2000, and put together a high-level and very comprehensive Editorial Board. Varadarajan and several other contributors to this volume served as associate editors from the very beginning and for many years, even after Chatterji retired. Professor Robert Dalang replaced him as the Managing Editor from 2001 onwards. Some may remember many fine articles that have appeared in his journal. A commentary on the mathematical diary of Gauss, translated and annotated, was one \cite{gray}. An article on the Poncelet porism \cite{bos} was another, and of course, there are many other examples of such articles, if one goes through the back issues.

One difficulty of running a journal that covers many different areas of mathematics is that the associate editors seldom meet at meetings, since most of them will not be present at the specialized meetings that each one often attends. To circumvent this, Chatterji would organize an editors meeting every four years at each International Congress of Mathematicians. Indeed, this would be one of the rare meetings that several editors would simultaneously attend. This editors meeting would be followed by a dinner at a nice local restaurant.

Against all predictions and misgivings of doomsayers, the journal has survived and is doing well. It publishes about 500 pages per year and is now an important source of expository articles.

\section{Service to the community}

Chatterji was vice-president of the Swiss Mathematical Society in 1984 and 1985, then president of this society in 1986 and 1987. He was president of Section VII of the Swiss Academy of Natural Sciences (this section covers mathematics, history of science and logic) from 1994 to 2000. From 1994 to 1998, Chatterji was a member of the ``Commission de Recherche de l'Acad\'emie suisse des sciences naturelles pour le Fonds National Suisse." He represented the Swiss Mathematical Society at meetings of the International Mathematical Union from 1986 to 2010.

Chatterji participated very actively in the organization of the 1994 International Congress of Mathematicians in Z\"urich. Indeed, he was Secretary of the Executive Committee for ICM 1994 (from 1991 to 1994), and he edited the two volumes of the proceedings of this congress \cite{rd14}.

Varadarajan recalls that  as a local organizer, Chatterji would be peppered with questions from people in all sorts of languages, and, much to everyone's surprise, he would reply in their own language almost every time! Some Italian colleagues even mentioned that when they spoke to him over the phone, they could not tell that he was not Italian!

Chatterji was a member of the Euler Committee of the Swiss Academy of Sciences from 1991 to 2000. This committee was originally established in 1907. Its goal, on which it has been working for over 100 years (!), was to bring to publication the entire scientific production of Euler, and moreover, to keep the published works {\em permanently available} (in 2010, the Euler-Gesellschaft was united with the corresponding one concerning the Bernoulli brothers and is now called the ``Bernoulli-Euler Center," at the University of Basel). This is a gigantic project, as anyone who is even remotely aware of how prolific Euler was, would understand. The committee has already published 80 volumes in quarto, and there are still four more volumes to come out, not to speak of his voluminous correspondence and other memorabilia. Here is what a co-member of the Euler Committee, Guido Mislin, wrote to Varadarajan:
\begin{quote}
Until my retirement from ETHZ in 2007, I was a member of the "Euler-Committee" and so was S.D.~Chatterji, for many years. His mathematical expertise and deep knowledge and understanding of history of mathematics made him a very valuable member. I remember vividly discussions on old and new mathematics that I had with  Chatterji and I always appreciated his wisdom and generosity. When we sat together at luncheons after committee meetings, I got to know him also from a non-mathematical perspective. I always was impressed by his knowledge on all kind of non-mathematical topics, his insight and his view on worldly matters. He was a great person.

We will all miss him.
\end{quote}
Varadarajan remembers an occasion when Chatterji talked about the Euler Committee. He said that a high level ministerial person asked him after one of the meetings, \lq\lq Professor Chatterji, could you tell me why we should be in this Euler business\rq\rq, to which Chat replied, \lq\lq Sir, if we do not do this, we will be the laughing stock of the whole mathematical world.\rq\rq This was Chat in a nutshell, never sugar-coating what should be done.

Chatterji was also a member of the board of the Stiftung zur F\"orderung der mathematischen Wissenschaften in der Schweiz (Foundation for the promotion of mathematical sciences in Switzerland), from 1991 to 2000. This foundation was created in 1929 by four mathematicians, including Michel Plancherel, with the objectives of encouraging young mathematicians, contributing to symposiums and supporting mathematical publications.

In addition to his editorial work with Expositiones Mathematicae, Chatterji was an associate editor of l'Enseignement Math\'ematique from 1989 to 2000. He was also the director of the collection ``Math\'ematiques" of the Presses Polytechniques et Universitaires Romandes from 1989 to 2000. This collection mainly published textbooks for EPFL students, as well as some higher level monographies.

   Finally, Chatterji also served as ``Commissaire" (expert) for various high schools in the canton of Ticino (1973 to 1989).

\section{The EPFL Mathematics library}

One feature that distinguished Chatterji from many other mathematicians was his exceptionally broad mathematical interests and culture. During his many early appointments, he met many different mathematicians and acquired substantial knowledge about most important mathematical subjects.  His colleagues would often go to him with questions about various topics, and he would usually supply a precise answer or a precise bibliographical reference. This broad interest may have been one reason why Chatterji became invested in the development of the mathematics library at EPFL.

This library had been founded in 1970 by Professor Bruno Zwahlen. In 1977, Chatterji became the head of the library, and from then on, had the final say in what books and journals were to be ordered by the library, and he remained in this position until his retirement in 2000. During this long reign, Chatterji made sure that the library would acquire each year all the important new math books. Under his direction, with the help of the dedicated library manager Madame Monique Gervaix who, from 1988 to 2007, supervised a small team of four librarians, the library became the leading math library in the French-speaking part of Switzerland with, in the year 2000, over 25,000 monographs and 12,500 bound volumes of periodicals, with about 500 regular users, of which more than half were from other departments and universities.

Chatterji felt passionately about his library. For instance, in 1982, when the math department moved from its downtown location in Lausanne to its current location in the suburban town of Ecublens, Chatterji securred for his library the best location in the new math building, on the third and top floor, with a magnificent view of the Lac L\'eman and the mountains behind it (even the snow-capped top of Mont Blanc was visible on a clear day). This was an inspiring setting, where professors and students felt ``at home" and enjoyed working and browsing. Books were classified using the AMS Subject Classification, making them particularly easy for users to find, and all books and journals were immediately accessible. When someone wanted to borrow a book, he would write his name on a small card that was fitted into the book and place this card in a designated box (this system was also in use at other universities). This system had two important advantages: if you wanted to consult a book that was checked out, you could immediately see who had it and go and ask the person for it. And you could see all the mathematicians who had previously checked out the book, and this could be quite a long and impressive list in the eyes of a young Ph.D.~student!

Chatterji was an unusually calm person, but his passion for the library would sometimes reveal a different facet of his personality. As recounted by Monique Gervaix, Chatterji had an unusual flair for detecting missing books, or books that were not placed in the correct place on the shelves, and this was one of the rare events that could really upset him! And he would also become upset when someone did not properly fill out the card for the book that they had borrowed (though only Madame Gervaix was aware of this)!

After his retirement, from 2001 to 2010, Chatterji became Director Emeritus of the library and continued to be the main person responsible for new acquisitions. He did this so well that when others would ask that the library order a new book, it would often turn out that the book was either already in the stacks or at least already ordered!

In 2010, the mathematics library was merged with all other EPFL libraries and the math collection was moved into a central location, though the classification and layout of the books was preserved. However, even after 2010, Chatterji continued to take care of acquisitions and did this up to the very last days of his life. Here are some words from Mr.~Julien Junod, who is the current math librarian, about the time when he worked with Chatterji:

\begin{quote}
 Prof. Chatterji's huge culture covered mathematics both in breadth (he felt comfortable in most fields of pure mathematics), and time (he was showing a great interest for older texts and the history of mathematics). This was a blessing for the library, whose acquisitions were chosen with great care and love. His choices were very rigorous and informed; they contributed to building a coherent and comprehensive book collection. The value of the library has been praised several times by international experts, and I think that it is one of the most beautiful math libraries of Western Europe.

 ...It was a great pleasure and opportunity to work with him. He was very friendly, calm, enthusiastic and interested - a very pleasant and reliable colleague! I learned a lot about mathematics, mathematicians, mathematical culture and history, just by looking at what books he selected, and of course also by talking to him. You could feel that he loved books deeply, and this showed in the great care and effort that he spent in developing the library...
\end{quote}

\section{The EPFL Mathematics Department}

Chatterji participated actively in the running of the EPFL Mathematics Department. We have just discussed his important work with the math library. He was also Chairman of the Department in 1976 and in 1994-95. He was a member of the Commission de Recherche of EPFL, which makes recommendations on research funding. This committee was in particular responsible for evaluating grant proposals of EPFL faculty members.

Chatterji participated actively in discussions about the mathematics curriculum for EPFL students, in particular for those majoring in mathematics. For many years, there were discussions among the faculty about how much time students should spend learning applied mathematics versus more theoretical mathematics, and Chatterji was one of the proponents of the theoretical side.

Chatterji was also a senior member of several hiring committees during the years 1990 to 2000. This was a new phase in which EPFL was aiming to attract high level international researchers. Chatterji participated in the appointments of many new mathematics professors, who would shape the EPFL mathematics department during the twenty ensuing years.

\section{Last years}

After he retired, Chatterji kept up his many interests in mathematics. He was a regular attendee of the S\'eminaire Bourbaki in Paris, France. He edited the collected papers of several mathematicians of the 19th century (\cite{rd16}--\cite{rd19}). He wrote on the work of, and controversies between, many nineteenth and twentieth century mathematicians \cite{c71,c72, c73b, u5,u6}. He became a well known personality among mathematical historians.
\begin{figure}[h]
\begin{center}
\includegraphics[height=7.5cm]{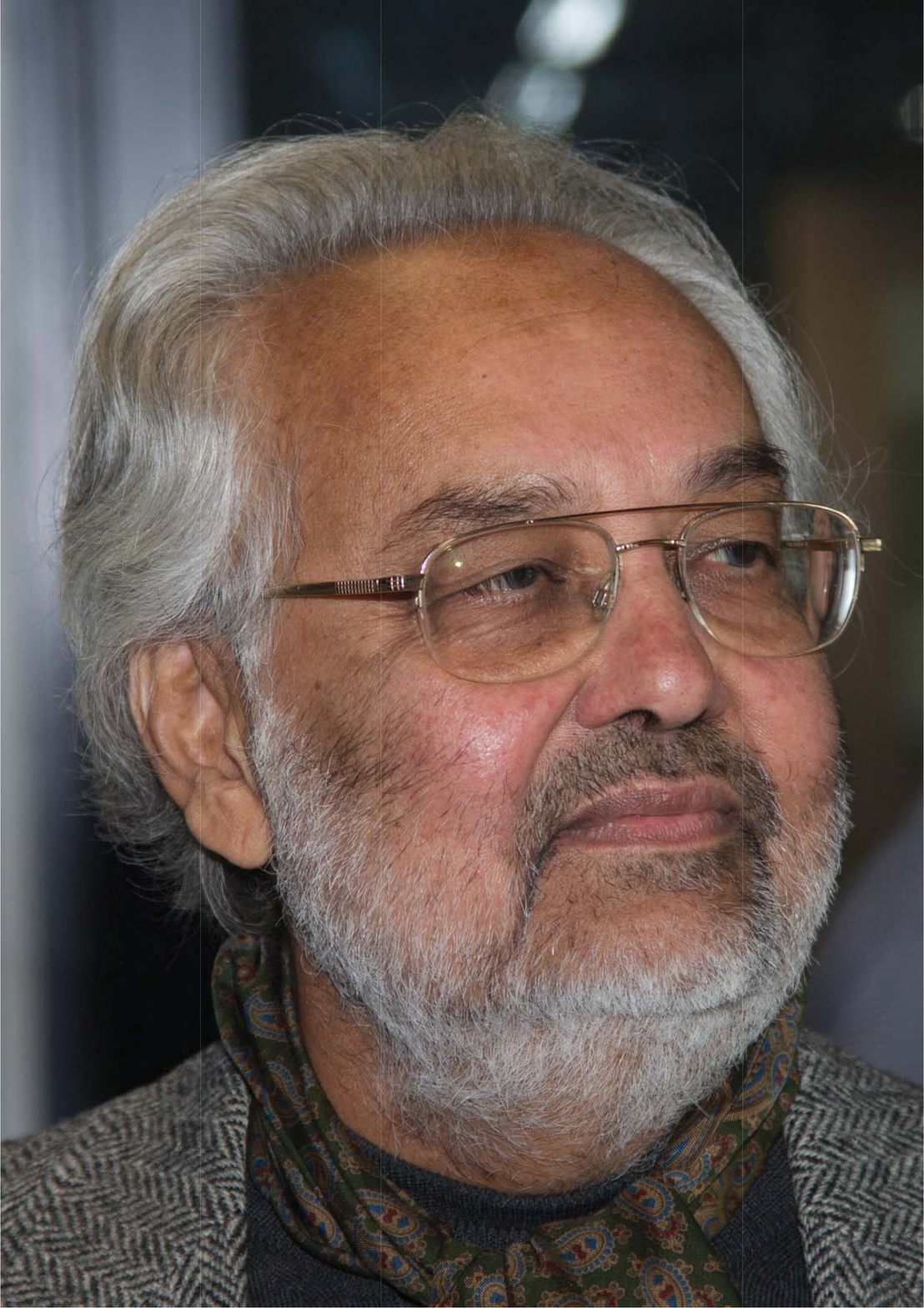}
\end{center}
\caption{Chatterji in 2010.  Photo by Prof.~Burchard Kaup (Universit\'e de Fribourg), printed with his permission.}
\end{figure}
Chatterji's main hobby was reading. His house was filled with books, both mathematical and nonmathematical. He felt that if people spent more time reading, then the world would be a better place.

He liked to watch certain English television series, and movies, such as those of the French filmmaker Jacques Tati. He was an avid football and tennis fan. He liked to listen to music, and had a fondness for Christmas carols, which reminded him of his youth. He enjoyed very much spending time in his wife's hometown near Locarno (Ticino).

In 2001, Chatterji underwent heart bypass surgery, as a precaution. The operation was successful and after a few months of convalescence, Chatterji soon returned to his usual schedule. This was followed by a second operation in 2008, which was also successful. In 2014, Chatterji suffered a mild stroke, which forced him to slow down and travel less. However, just a few weeks before his death, Chat and Carla went to Italy to celebrate their fiftieth wedding anniversary!

On Monday, September 25, 2017, Chatterji went into the library as usual and worked with Julien Junod. Tuesday morning, he had difficulty breathing, and was taken to the main hospital in Lausanne, where he passed away two days later, due to a combination of respiratory and cardiac problems. Chatterji is survived by his wife Carla, his daughter Indira and two grandchildren.

\section{Some personal comments of Varadarajan}

Chatterji and I spoke with each other about every other month. Over these long phone conversations, we explored all topics. He was of course more erudite than I, but that was never an obstacle. His father was long-lived, and so we were certain he would outlive me to write my obituary, and so he insisted that every piece of work that I published should be sent to him. I faithfully kept my end of the bargain, but why should it be that he had to go first?

As yet another example of his diversity, kindness, and generosity, I may say something personal. One of my grand nieces lives in Switzerland. She has aspirations for a literary career, and has already published many poems in English as well as translations of Urdu poetry. When the time came for her to continue her studies, I suggested she meet with Chat. He was very generous with his time. He found her achievements impressive, and gave her good advice. Here is what she wrote to me after I told her about Chat's passing away:
\begin{quote}
I am so sorry to hear that. Our limited interaction was wonderful and he was such a warm and well-spoken man. I will remember him fondly.
\end{quote}

\noindent In a note to me, my friend Joop Kolk said the following:
\begin{quote}
My memories of meeting with Chatterji, in Lausanne and Zurich, are of a sophisticated man, highly cultured and very knowledgeable, but also of a gentle and generous disposition. I had the impression that he and I quickly liked each other. Later on we had a most agreeable interaction by e-mail on some mathematical issues. I can very well understand that losing such a friend is most painful for you.
\end{quote}

\noindent And here is an excerpt from a letter that K.R.~Parthasarathy, wrote to me after learning that Chat was no more:
\begin{quote}
 I came to know his name through his martingale limit theorems in Banach spaces. Later in life I met him in Lausanne and Zurich. At the 1994 ICM in Zurich he took personal care of me because he was very much aware of the difficulties an Indian of my kind would experience in a country like Switzerland. He had visited ISI, Delhi, given seminars and visited our flat in the campus with Swiss chocolates for our boys. He was fluent in many European languages and deeply interested in philosophy. He loved good wine. Now he is a page in the book of my cherished memories. This book too is aging. His daughter Indira had also given seminars at ISI as well as the neighboring Jawaharlal Nehru University.
\end{quote}

Chat lived his life fully, enjoying all aspects. He leaves behind an unfinished autobiography (which, he told me, would be entitled {\em Autobiography of an unknown Indian mathematician}, in conscious imitation of Nirad Chowdury's famous book {\it Autobiography of an unknown Indian}). Of course he was not all that unknown as the title presumes. He has also left behind an unfinished manuscript on Measure Theory. Given his knowledge of {\it all} aspects of Measure Theory and its applications, this would have been {\it the} book on Measure theory.

He was a great teacher, not only at the basic level, but also as a supervisor to doctoral students. Professsor Francesco Russo, who was a student of Chatterji, wrote to me that Chatterji was like a father to him. He was the ideal illustration of the phrase {\it in loco parentis}: see the testimonial article by F.~Russo in this volume.

Chatterji was in a deep sense a true renaissance man. He loved literature, especially poetry. He would recite the opening lines of Milton's {\em Paradise Lost} with a feeling that would thrill the listener to the very core. I cannot describe how many hundreds of hours we talked, about every topic on earth, and how each time we talked, the occasion would bring inconceivable pleasure to me.

To be born is to die, and Chatterji's life was a remarkable illustration of the philosophy that the journey is the thing, and that being fully alive is the extraordinary phenomenon. The quotation at the beginning of this article is from Ingmar Bergman's film {\it Seventh seal}, which celebrates life amid the devastation of death caused by the Great Plague in medieval Sweden. It is an apt description of his entire life, and the statement of J\"ons, the squire of the knight in that film, is one that Chatterji could very well have made; in fact did make indirectly by the exemplary life he led, as a husband, father, creative thinker, scholar, teacher, friend, philosopher, and guide.
\vskip 16pt

\noindent{\sc Acknowledgement.} This article was initially based on a preprint prepared by V.S.~Varadarajan. Additional material was added by R.C.~Dalang, using information provided to R.C.D.~by Chatterji himself upon his retirement in 2000. R.C.D.~also consulted Chatterji's annual scientific reports that had been preserved by his family, and included information supplied by Carla Chatterji-Bolognini and Indira Chatterji, by Shyam Sunder Chatterjee, Dharani Dhar Chatterjee and Dipika Chatterjee, as well as by Chatterji's son-in-law Fran\c{c}ois Labourie. Librarians Julien Junod and J\'er\^ome Yerly helped identify certain of Chatterji's publications, and Professor Burchard Kaup provided several photographs of Chat. We thank all of these kind people for their cooperation.

\eject

{\small
\def\refname{References}

\noindent{\bf Ph.D.'s supervised by S.D.~Chatterji}
\vskip 16pt

\begin{itemize}

\item[] Lischer, Peter. Robust regression.  Th\`ese no.~160. \'Ecole Polytechnique F\'ed\'erale de Lausanne (1973)

\item[] Brunnschweiler, Andreas.  Mathematische Untersuchung eines Diffusionsmodells; die L\"osung der Gleichung von Ornstein-Uhlenbeck als Limes von L\"osungen der Rayleigh Gleichung. Th\`ese no.~260. \'Ecole Polytechnique F\'ed\'erale de Lausanne (1976)

\item[] Carnal, Etienne. Processus markoviens \`a plusieurs param\`etres. Th\`ese no.~325. \'Ecole Polytechnique F\'ed\'erale de Lausanne (1979)

\item[] Russo, Francesco. Champs markoviens et pr\'ediction. Th\`ese no.~707. \'Ecole Polytechnique F\'ed\'erale de Lausanne (1987)

\item[] Bolle, Yves. Int\'egrales de Riemann dans les espaces abstraits. Th\`ese no.~2107. \'Ecole Polytechnique F\'ed\'erale de Lausanne (2000)

\item[] Abaza, Riadh. Int\'egrale de Riemann et continuit\'e dans les espaces vectoriels topologiques. Th\`ese no.~2147. \'Ecole Polytechnique F\'ed\'erale de Lausanne (2000)

\end{itemize}

\def\refname{Publications of S.D.~Chatterji\footnote{This list of publications was prepared using five sources: a list compiled by Varadarajan, a list prepared by Chatterji himself, but which stopped in 1998 and only contained his research papers, MathSciNet, Google Scholar, and the EPFL scientific data base Infoscience. None of these sources seem to list all the documents mentioned here.}}

}
\vskip 16pt

\begin{tabular}{ll}
V.S.~Varadarajan & Robert C.~Dalang \\
Department of Mathematics & Institut de math\'ematiques \\
University of California Los Angeles& \'Ecole Polytechnique F\'ed\'erale de Lausanne \\
Los Angeles, CA 90024& Station 8 \\
U.S.A.& CH-1015 Lausanne \\
vsv@math.ucla.edu & Switzerland \\
 & robert.dalang@epfl.ch
\end{tabular}

\end{document}